\documentclass[12pt]{iopart}

\usepackage{iopams}
\usepackage{graphicx}
\usepackage{soul}
\usepackage{multirow}
\usepackage{booktabs}
\usepackage{diagbox}
\usepackage{makecell}
\usepackage{filecontents}
\usepackage[dvipsnames]{xcolor}
\usepackage{hyperref}
\usepackage{bm}

\newcommand\myshade{85}
\colorlet{mylinkcolor}{BrickRed}
\colorlet{mycitecolor}{NavyBlue}
\colorlet{myurlcolor}{Aquamarine}

\hypersetup{
  linkcolor  = mylinkcolor!\myshade!black,
  citecolor  = mycitecolor!\myshade!black,
  urlcolor   = myurlcolor!\myshade!black,
  colorlinks = true,
}

\newcommand{\Mod}[1]{\ (\mathrm{mod}\ #1)}

\expandafter\let\csname equation*\endcsname\relax
\expandafter\let\csname endequation*\endcsname\relax
\usepackage{amsmath}
\usepackage{cleveref}

\newtheorem{definition}{Definition}
\newtheorem{example}{Example}

\DeclareMathOperator{\TC}{TC}
\DeclareMathOperator{\DTC}{DTC}
\newcommand{\bX}{\bi{X}}
\newcommand{\qed}{\hfill$\blacksquare$}

\newcommand{\rev}[1]{\textcolor{black}{#1}} 

\makeatletter
\long\def\@makefntext#1{\parindent 1em\noindent 
 \makebox[1em][l]{\footnotesize\rm$\m@th{^\arabic{footnote}}$}%
 \footnotesize\rm #1}
\def\@makefnmark{\footnotesize{${^\arabic{footnote}}\m@th$}}
\def\@thefnmark{\arabic{footnote}}
\makeatother
\usepackage{natbib}
\setcitestyle{authoryear,open={(},close={)}}

\begin{document}	
	\title[]{Hyperharmonic analysis for the study of high-order information-theoretic signals}
	
	\author{Anibal M. Medina-Mardones$^{1,2}$, Fernando E. Rosas$^{3,4,5}$, Sebasti\'{a}n E. Rodríguez$^{6}$, Rodrigo Cofr\'{e}$^{7}$}
	
	\address{
	$^{1}$ \quad Laboratory of Topology and Neuroscience, \'{E}cole Polytechnique F\'{e}d\'{e}rale de Lausanne, Lausanne, Switzerland \\
	{$^{2}$} \quad Department of Mathematics, University of Notre Dame du Lac, Notre Dame, Indiana, USA \\
	$^{3}$ \quad Data Science Institute, Imperial College London, London SW7 2AZ, UK \\
	$^{4}$ \quad Centre for Psychedelic Research, Department of Medicine, Imperial College London, London SW7 2DD, UK \\
	$^{5}$ \quad Centre for Complexity Science, Imperial College London, London SW7 2AZ, UK \\
	$^{6}$ \quad Universidad T\'{e}cnica Federico Santa Mar\'{i}a, Departamento de Inform\'{a}tica, Valpara\'{i}so, Chile\\
	$^{7}$ \quad CIMFAV-Ingemat, Facultad de Ingenier\'{i}a, Universidad de Valpara\'{i}so, Valpara\'{i}so, Chile
	}
	\ead{anibal.medinamardones@epfl.ch}
	\vspace{10pt}
	\begin{indented}
		\item[]February 2021
	\end{indented}

\begin{abstract}
Network representations often cannot \rev{fully} account for
the structural richness of complex systems \rev{spanning}
multiple levels of organisation. Recently proposed high-order information-theoretic signals are well-suited to capture synergistic phenomena that transcend pairwise interactions; however, the exponential-growth of their cardinality severely hinders their applicability. In this work, we combine methods from harmonic analysis and combinatorial topology to construct efficient representations of high-order information-theoretic signals. 
The core of our method is the \rev{diagonalisation} of a discrete version of the Laplace-de Rham operator, that geometrically encodes structural properties of the system. We capitalise \rev{on} these ideas by developing a complete workflow for the construction of hyperharmonic representations of high-order signals, which is applicable to a wide range of scenarios.
\end{abstract}

\vspace{4pc}
\noindent{\it Keywords}: high-order phenomena, Laplace operator,
harmonic analysis, signal processing, information theory
	
\newpage
\section{Introduction}
	
The principle of representing interdependencies as networks \rev{has} revolutionised complexity science by introducing a systematic approach to gain insight \rev{into} 
the inner structure of a wide range of complex systems~\citep{newman2018networks}.
These networks provided a \textit{lingua franca} to describe the properties of interdependencies found in chemical, biological, social, and technological systems~\citep{vasiliauskaite2020understanding}, and have enabled great advances in a widening range of areas including computational neuroscience~\citep{Rubinov2010}, human evolution~\citep{donges2011nonlinear}, financial analysis~\citep{Bonanno2004}, and epidemic spreading~\citep{Pastor-Satorras2001}, just to name a few.
However, by their very nature, these methods focus on the analysis of pairwise interactions, and hence are prone to \rev{miss} 
important high-order synergistic phenomena that are \rev{a hallmark of many}
complex systems. 

This critical limitation of traditional network analyses has been acknowledged by a number of recent research efforts that focus on developing techniques to study \rev{high}-order
interactions~\citep{Petri2014,Iacopini2019,Battiston2020}.
These developments have provided novel 
techniques capable of, for example, detect\rev{ing} non-local structures~\citep{Petri2013}, highlight\rev{ing} the role of inhomogeneities in functional connections~\citep{Petri2014}, and character\rev{ising} discontinuous transitions ~\citep{Iacopini2019}. 
However, it is important to notice that many of these approaches are based on hypergraphs and other high-order structures \rev{built} solely from pairwise statistics, and hence their scope remains limited.

An attractive set of tools to further extend the reach of high-order analyses can be found in a parallel body of work, which originated in efforts to develop tools to capture high-order statistical phenomena related to the brain~\citep{tononi1994measure,schneidman2003synergy,latham2005synergy,ganmor2011sparse}. Particularly interesting are multivariate extensions of Shannon's mutual information, including the Interaction Information~\citep{mcgill1954multivariate}, Total Correlation \citep{watanabe1960information}, and Dual Total Correlation \citep{Han1978}, which can be used to gain insights about the high-order structure \rev{exhibited by} 
groups of three or more interdependent variables (see e.g.~\cite{timme2014synergy,baudot2019topological}). In this work we focus on the recently proposed \textit{\mbox{O-information}}~\citep{rosas2019quantifying}, which is a principled tool to identify synergy-dominated systems, and has \rev{been} found to be relevant for analysing various complex systems --- including the study of neural spiking data~\citep{stramaglia2020quantifying} and aging in fMRI data~\citep{gatica2020high}. 

When applied to large systems, metrics such as the O-information are naturally represented as high-order signals whose domain is the set of all hyper-edges on a regular hypergraph.
Because the cardinality of these hypegraphs grows super-exponentially with the system size, a key open problem is to find efficient ways to represent the content of these signals. While simple approaches, such as computing the average of the signal across dimensions can be effective 
\citep{gatica2020high}, an important challenge is to find principled ways to generate low-dimensional representations of these signals while preserving their intrinsic relational structure across the hypergraph. 
A popular technique to study similar issues in the case of traditional (weighted) networks is to transform the signals into a basis of eigenvectors of the graph Laplacian \citep{belkin2003laplacian}, which has been used with great success in the context of graph signal processing (see e.g.~\cite{shuman2013emerging, sandryhaila2014discrete, Atasoy2016, atasoy2017connectome, expert2017graph}). Related methods based on harmonic analysis and combinatorial topology, which we refer to as \textit{hyperharmonic analysis}, have been recently developed to study high-order signals (see e.g.~\cite{barbarossa2020topological}); however, they have not yet --- to the best of our knowledge --- been used to analyse high-order information-theoretic signals.
Other aspects of the relationship between high-order information-theoretic quantities and topology \rev{have} been explored in~\cite{baudot2015homological,baudot2019poincare}.

This article establishes a bridge between the domains of high-order information-theory and hyperharmonic analysis, and introduces a complete workflow for the study of high-order information-theoretic signals using the hyperharmonic modes of a structural simplicial complex.
Our choice of hyperharmonic modes is based on a discrete version of the Laplace-de Rham operator \rev{that }
geometrically encodes the strength of low-order interactions. 
Our workflow makes no assumptions about the structure of the data, and hence can be applied to a broad range of scenarios. As a proof of concept, we illustrate our approach \rev{by} analysing the musical scores of the latter symphonies written by F.J. Hadyn, \rev{where} our results demonstrate the far superior dimensionality-reduction capabilities of our method compared to other representations.

The rest of the article is structured as follows. Section~\ref{sec:hoitm} introduces key notions from multivariate information theory, reviewing the state-of-the-art of high-order metrics.
Then, Section~\ref{sec:hyper} presents 
fundamental notions from combinatorial topology required to define the Fourier transform of \rev{a} higher-dimensional signal. Section~\ref{sec:pipe} introduces our proposed workflow, and Section~\ref{sec:example} illustrates the method \rev{used} on Hadyn's symphonies. Finally, Section~\ref{sec:conc} summarizes our conclusions and discusses future work.
	
\section{High order information-theoretic measures} \label{sec:hoitm}
	
Let us consider a scientist studying a complex system, whose state is described by the vector $\bi{X}^N = (X_0, \dots, X_N)$. 
Let us assume that the scientist has enough data to allow for the construction of a reliable statistical description of its joint probability, which is denoted by $p(X_0, \dots, X_N)$. The goal of this study is to leverage the statistics encoded in $p$ in order to understand the structure of interdependencies that characterize $\bi{X}^N$. This endeavor can lead either to \rev{the building of} statistical markers to classify different systems --- or different states of the same system, or to \rev{the building of} parallels between seemingly heterogeneous systems based on the similarity of their internal structure.

Through this section, random variables are denoted by capital letters (e.g. $X,Y$) and their realisations by lower case letters (e.g. $x,y$). Random vectors and their realisations are denoted by capital and lower case boldface letters, respectively.

\subsection{\rev{Networks based on pairwise statistics}}
\label{sec:pairwise}

A popular way to analyze the interactions within $\bi{X}^N$ is to represent them as networks, where each variable $X_0, \dots, X_N$ is represented as a node, and edges between nodes represent the strength of their interaction. A simple way to build such a network is to calculate the correlation matrix $\mathcal{R}_{\bi{X}^N} := [r_{i,j}]$ with components given by
\begin{equation}
 r_{i,j}=  \sum_{x_i,x_j} p(x_i,x_j) \Big[x_i x_j - \sum_{x_i}x_i \sum_{x_j}x_j  \,\Big] .
\end{equation} 
and use it as an adjacency matrix --- either binarising its components via a threshold, or considering weighted edges. This construction, however, only captures linear \rev{dependencies} between the variables. 
More encompassing analyses often consider non-linear measures of dependency such as Shannon's mutual information and focus on the matrix $\mathcal{I}_{\bi{X}^N} := \big[ I(X_i;X_j) \big]$ of mutual information terms, which are computed as
\begin{equation}
    I(X_i;X_j) = \sum_{x_i,x_j} p(x_i,x_j) \log  \frac{p(x_i,x_j)}{p(x_i) p(x_j)}.
\end{equation}
\rev{However, please note that this description still does not fully characterises the interdependencies of a multivariate system. In particular, as discussed in the next section, a system can have zero pairwise dependencies (i.e. $I(X_i;X_j)=0$ for $i,j\in\{1,\dots,N\}$ for all $i\neq j$) while still having high-order interdependencies.}

\subsection{High-order statistical effects}
\label{sec:high-order}

It is important to realise that the joint probability distribution $p(\bi{X}^N)$ may contain substantial information that is not assessed by any of the pairwise marginals $p(X_i,X_j)$. An elegant way of gaining insights into this is provided by information geometry, as presented in \cite{Amari2001} (for related discussions, see also~\cite{rosas2016understanding}). Consider the $k$-marginals that are obtained by marginalising $p(\bi{X}^N)$ over $N-k$ variables. One can see that the $k$-marginals provide a more detailed description of the system than the $(k-1)$-marginals by noting that the latter can be directly computed from the former by marginalising the corresponding variables. In contrast, the process of marginalising involves irreversible information loss, as there are many $k-1$ marginals that are consistent with a given  set of $k$-marginals. 

Perhaps the simplest example of high-order statistical dependency is given by the exclusive-or (xor) logic gate. Consider two independent fair coins, $X_0$ and $X_1$, and let 
\begin{equation}
X_2 = X_0 \texttt{(xor)} X_1 =     
\begin{cases}
0 \quad \text{if} \,\,\, X_0 = X_1, \\
1 \quad \text{otherwise.}
\end{cases}
\end{equation} 
A quick calculation shows that the mutual information matrix $\mathcal{I}_{\bi{X}^2}$ of $(X_0, X_1, X_2)$ has all its off-diagonal elements equal to zero, making it indistinguishable from an alternative situation where $X_2$ is just another independent fair coin. Therefore, put simply, any ``\mbox{\texttt{xor}-like}'' effects \rev{are} completely ignored by constructions based only on pairwise statistics --- either networks or hypergraphs.

While commonly neglected, high-order statistics \rev{have recently been proven} to be instrumental in a number of systems. For example, \rev{high-order phenomena have been shown to play a key role in self-organising} distributed systems, including the most complex types of elementary cellular automata~\citep{rosas2018information}.  \rev{It has also been suggested that high-order interactions could drive thermalisation processes within closed systems~\citep{lindgren2017approach}.} Additionally, high-order interactions \rev{have} been argued to play a key role in neural information processing~\citep{wibral2017partial} and high-order brain functions~\citep{tononi1994measure}, being at the core of popular metrics employed in consciousness science~\citep{mediano2019beyond}. Furthermore, it has recently been shown that high-order statistics also play a crucial role in enabling emergent phenomena~\citep{rosas2020reconciling}.

\subsection{O-information}
\label{sec:Oinfo}

Shannon's mutual information is limited to \rev{capturing} the dependencies of two groups of variables, but cannot directly assess triple or higher interactions. The most popular non-negative multivariate extensions of the mutual information are the \textit{Total Correlation} (TC)~\citep{watanabe1960information} and the \textit{Dual Total Correlation} (DTC)~\citep{Han1978}, which are defined as
\begin{align*}
\TC(\bX^N) &:= \sum_{i=0}^N H(X_i) - H(\bX^N) , \\
\DTC(\bX^N) &:= H(\bX^N) - \sum_{i=0}^N H( X_i \mid \bX^N_{-i}) .
\end{align*}
Above, $H(X_i) =-\sum_{x_i} p(x_i) \log p(x_i)$ corresponds to Shannon's entropy, $H(X_i|X_j) = H(X_i,X_j) - H(X_j) $ is the conditional Shannon entropy, 
and $\bX^N_{-i}$ is the vector of all variables except $X_i$ (i.e., $\bX^N_{-i}=\left( X_0, \ldots, X_{i-1}, X_{i+1}, \ldots, X_N \right)$). Importantly, both $\text{TC}$ and $\text{DTC}$ are zero if and only if all variables $X_0, \dots, X_N$ are jointly statistically independent --- i.e. if $p(\bi X^N) = \prod_{i=0}^N p(X_i)$.

Unfortunately, both $\text{TC}$ and $\text{DTC}$ provide 
metrics for high-order interdependency which are difficult to analyse together. A recent approach, introduced in~\cite{rosas2019quantifying}, proposes to employ a linear transform over these two metrics to obtain the \textit{O-information} and the \textit{S-information}, which have more intuitive interpretations.

\begin{definition}  Given a set of $N+1$ random variables $\bX^N = (X_0, \dots, X_N)$, their \mbox{O-information} is defined as
	\begin{equation}\label{oinfo}
	\Omega(\bX^N) = \TC(\bX^N) - \DTC(\bX^N).
	\end{equation}
Similarly, their S-information is defined as
	\begin{equation}\label{sinfo}
	\Sigma(\bX^N) = \TC(\bX^N) + \DTC(\bX^N).
	\end{equation}
\end{definition}

The O-information can be seen as a revision of the measure of neural complexity proposed by Tononi, Sporns and Edelman in~\cite{tononi1994measure}, which provides a mathematical construction that is closer to their original desiderata~\citep{rosas2019quantifying}. \rev{In effect, the O-information is a signed metric that captures the balance between high- and low-order statistical constraints. While low-order constraints impose strong restrictions on the system and allow little amount of shared information between variables, high-order constraints impose collective restrictions that enable large amounts of shared randomness.\footnote{For a more technical discussion, see Sec.2 in ~\citep{rosas2019quantifying}).} In particular, high-order constraints can generate global interdependencies that do not impose corresponding pairwise dependencies, as observed e.g. in the \texttt{xor} logic gate studied in Section~\ref{sec:high-order}.}

\rev{By construction, $\Omega (\bi X^N) < 0$ implies a predominance of high-order constraints within the system $\bi X^N$, a condition that is usually referred to as \textit{statistical synergy}. Conversely, $\Omega (\bi X^N) > 0$ implies that the system $\bi X^N$ is dominated by low-order constraints, which imply \textit{redundancy} of information. This nomenclature is further supported by the following key properties of the O-information}:
\begin{itemize}
    \item[(1)] \rev{It captures genuine high-order effects, at it is zero for systems with only pairwise interdependencies: if the joint distribution of $\bX^{N-1}$ (for $N$ odd) can be factorised as $p_{\bX^N}(\boldsymbol x^N) = \prod_{k=0}^{(N-1)/2}p_{X_{2k},X_{2k+1}}(x_{2k},x_{2k+1})$, then $\Omega(\bX^N)=0$.}
    \item[(2)] \rev{The O-information is maximised by redundant distributions where the same information is copied in multiple variables, and is minimised by synergistic (``\mbox{\texttt{xor}-like}'') distributions: e.g. for binary variables, $\Omega$ is maximised by the ``n-bit copy'' where $X_1$ is a Bernoulli r.v. with parameter $p=1/2$ and $X_0=\dots=X_{N-1}$, and is minimised when $X_0,\dots,X_{N-2}$ are i.i.d. fair coins and $X_{N-1}=\sum_{j=0}^{N-2} X_j \Mod{2}$.}
    \item[(3)] \rev{The O-information characterises the dominant tendency, being additive over non-interactive subsystems: if the system can be factorised as $p_{\bX^N}(\boldsymbol x^N) = p_{X_1,\dots,X_{m}}(x_1,\dots, x_m) \times p_{X_{m+1},\dots,X_N}(x_{m+1},\dots, x_N)$, then $\Omega(\bX^N)= \Omega(X_0,\dots,X_{m-1}) + \Omega(X_m,\dots,X_{N-1})$.}
\end{itemize}
\rev{Further insights, properties, and mathematical proofs related to the O-information can be found in~\citep{rosas2019quantifying}.}

On the other hand, the S-information is an over-encompassing account of interdependencies taking place at all orders --- being sometimes described as  ``very mutual information'' \citep{james2011anatomy}. In fact, a quick calculation shows that the S-information can be decomposed as $\Sigma(\bi X^N) = \sum_{i=0}^{N-1} I(X_i; \bm X_{-i}^N)$, which can be seen as a chain rule where the interdependencies involving each variable are sequentially addressed. \rev{By comparing this expression with the description of interdependencies given by pairwise mutual information terms discussed in Section~\ref{sec:pairwise}, one can see $S$ as similar to a summation of the links' strengths in standard network analysis --- however, in this case actually including the weight corresponding to beyond-pairwise interdependencies. Therefore, $\Sigma$ provides a scalar metric that assesses the overall strength of the interdependencies within a system; in particular, in contrast with pairwise metrics (c.f. Section~\ref{sec:pairwise}), it can be shown that $\Sigma(\bX^N)=0$ if and only if the parts of a system $\bX^N$ are statistically independent.}

In summary, while the $\text{TC}$ and $\text{DTC}$ provide alternative representations to the same construct, $\Omega$ and $\Sigma$ provide a complementary account of the system:
the latter address\rev{ing} the overall strength of interdependencies, and the former qualitatively characteris\rev{ing} the\rev{ir} dominant nature.

\subsection{Other metrics of high-order effects} 
\label{subsec: other high order signals}

Another popular metric of high-order interdependencies \rev{is} the \textit{Interaction Information}, first introduced in~\citep{mcgill1954multivariate} for systems with three variables.\footnote{The Interaction Information is closely related to the \textit{I-measures}~\citep{yeung1991}, the \textit{co-information}~\citep{Bell2003}, and the {multi-scale complexity}~\citep{Bar-Yam2004}} Building on an application of the inclusion-exclusion principle to entropies, the Interaction Information of $\bi X^n$ is a signed metric given by

\begin{equation}\label{eq:interaction}
I(X_0;X_1;\dots;X_N) := - \sum_{\bi{\gamma}\subseteq\{0,\dots,N\}} (-1)^{|\bi{\gamma}|} H(\bi{X}^{\bi{\gamma}}).
\end{equation}
where the sum is performed over all subsets $\bi{\gamma} \subseteq \{0, \dots, N\}$, with $|\bi{\gamma}|$ the cardinality of $\bi{\gamma}$, and $\bi{X}^{\bi{\gamma}}$ the vector of all variables with indices in $\bi{\gamma}$. While this measure has a direct interpretation as redundancy minus synergy for $N=2$, it no longer reflects this balance for larger system sizes~\cite[Section V]{williams2010nonnegative}. However, the Interaction Information can still be interpreted for arbitrary $N$ under a topological formulation of information, as described in~\cite{baudot2015homological,baudot2019topological}.

Other well-known metrics of high-order effects include the \textit{Redundancy-Synergy Index}~\citep{chechik2002group,timme2014synergy}, the \textit{Connected Information}~\citep{Amari2001,schneidman2003network}, and the \textit{Partial Entropy Decomposition}~\citep{ince2017partial}. Furthermore, a detailed exploration of multivariate decompositions can be found in the \textit{Partial Information Decomposition} framework~\citep{williams2010nonnegative}
and its constantly growing associated literature~\citep{faes2017multiscale,finn2018pointwise,ay2019information,rosas2020operational,makkeh2020differentiable}. 

\rev{Finally, it is worth mentioning \textit{Transfer Entropy} (TE)~\citep{schreiber2000measuring,barnett2009granger}, the non-linear generalisation of the well-known Granger causality\citep{granger1969investigating}. TE is a popular method to build directed network representations from time series, which has shown great effectiveness for discriminating features of interest in practical scenarios (see e.g.~\cite{marinazzo2014information,seth2015granger,bossomaier2016introduction,deco2021revisiting}). However, recent investigations have been noted that TE conflates different types of high-order phenomena~\citep{james2016information,mediano2019beyond}, making it usage for high-order analyses not straightforward. Interestingly, the tools introduced by~\cite{mediano2019beyond} allow to disentangle these heterogeneous effects, whose exploration constitutes an active area of research.}

\section{Hyperharmonic analysis}\label{sec:hyper}

This section, subdivided into three parts, introduces the basic concepts from combinatorial topology that we used to decompose high-order signals into hyperharmonic modes.
First, Section~\ref{sec:31} describes the objects over which signals will be considered; these are higher-dimensional versions of weighted graphs known as weighted simplicial complexes. Then, Section~\ref{sec:32} introduces the algebraic structure used to model high-order signals on weighted simplicial complexes; it consists of a family of inner product spaces, one for each dimension, and a pair of canonical linear maps between adjacent spaces. Finally, Section~\ref{sec:33} introduces the discrete analogue of the Laplace-de Rham operator, and defines the Fourier basis using a maximal set of linearly independent eigenvectors of this operator. Throughout the presentation, we provide references to more general treatments and original sources when possible.

\subsection{Simplicial complexes}
\label{sec:31}

A hypergraph $S = (V, E)$ is determined by a set of vertices $V = \{0, \dots, N\}$ with $N\in\mathbb{N}$ and a set of hyper-edges $E\subseteq \mathcal{P}(V)$, where $\mathcal{P}(V)$ is the power set of $V$. Furthermore, a hypergraph $S$ is said to be a \textit{simplicial complex} if it satisfies two conditions:
\begin{itemize}
    \item[(1)] all singletons $\{k\}$ with $k\in V$ are included in $E$, and
    \item[(2)] if $\sigma \in E$ and $\rho$ is a subset of $\sigma$, then $\rho \in E$.
\end{itemize}
Please note that the passage to simplicial complexes does not restrict the theory of hypergraphs significantly, since to every hypergraph one can assign a canonical simplicial complex by \textit{downward closure}. Explicitly, this is the smallest simplicial complex that contains the hypergraph. 
For an introduction to graphs and hypergraphs, we refer to \cite{berge1973graphs}; for a comprehensive introduction to hypergraphs in the context of complex system analysis see~\cite{johnson2013hypernetworks}.

If $S = (V, E)$ is a simplicial complex, the elements of $E$ are referred to as \textit{simplices}; and their \textit{dimension} is defined as one less than their cardinality (i.e. the number of
\rev{elements} they connect). This shift can be understood noting that the natural dimension of a point, corresponding to a singleton, is $0$.
The subset of $E$ consisting of simplices of dimension $n$ is denoted by $S_n$. The elements of $S_n$ are typically called ``$n$-simplices'', with the $0$-simplices and $1$-simplices being informally called vertices and edges, respectively.
If $v_0 < \dots < v_n$ are the vertices of a simplex, we denote this simplex by $\lbrack v_0, \dots, v_n \rbrack$.
\begin{example} \label{example: standard simplex}
The simplicial complex $\Delta^N$ with vertices $V = \{0, \dots, N\}$ and containing all possible simplices is referred to as the \textit{standard $N$-simplex}. In our applications, the $N+1$ vertices of the $\Delta^N$ will be associated with a set of random variables $X_0, \dots, X_N$. \qed
\end{example}

\subsection{Higher-dimensional signals}
\label{sec:32}

\rev{Following the terminology of \citep{barbarossa2020topological} and \citep{schaub2018flow}}, an $n$-\textit{dimensional signal} on a simplicial complex $S$ is \rev{defined to be} a function $\alpha : S_n \to \mathbb R$, that is to say, an assignment of a real number to each $n$-simplex. We refer to $n$-signals with $n \geq 2$ as \textit{high-order signals}.

\begin{example} \label{example: omega and sigma as signals}
	Consider a collection of random variables $X_0, \dots, X_N$. For every $n \in \{2, \dots, N\}$, the high-order $n$-signals $\Omega$ and $\Sigma$ on $\Delta^N$ are defined by
	\begin{align*}
	\Omega([v_0, \dots, v_n]) & := \Omega(X_{v_0}, \dots, X_{v_n}), \\
	\Sigma([v_0, \dots, v_n]) & := \Sigma(X_{v_0}, \dots, X_{v_n}).
	\end{align*}
	\qed
\end{example}

For a simplicial complex $S$ and $n \in \mathbb{N}$, we denote by $C_n(S)$ the vector space generated by all the $n$-simplices of $S$, i.e., \rev{the vector space whose elements are formal linear combinations of simplices}
\begin{equation*}
C_n(S) = \{\alpha_1 \sigma_1 + \cdots + \alpha_r \sigma_r\ |\ \alpha_i \in \mathbb R \text{ and } \sigma_i \in S_n\}.
\end{equation*}
Please note that there is a natural bijection between the set of $n$-signals on $S$ and $C_n(S)$ established by
\begin{equation*}
\alpha \mapsto \sum_{\sigma \in S_n} \alpha(\sigma)\, \sigma.
\end{equation*}
\rev{T}his bijection is used implicitly when referring to the elements of $C_n(S)$ as $n$-signals.\footnote{In the mathematics literature,
\rev{the elements of $C_n(S)$} are
called ``real-valued \rev{$n$}-chains'', but we do not use this terminology.}

\begin{figure}[h!]
    \centering
    \includegraphics[scale=0.2]{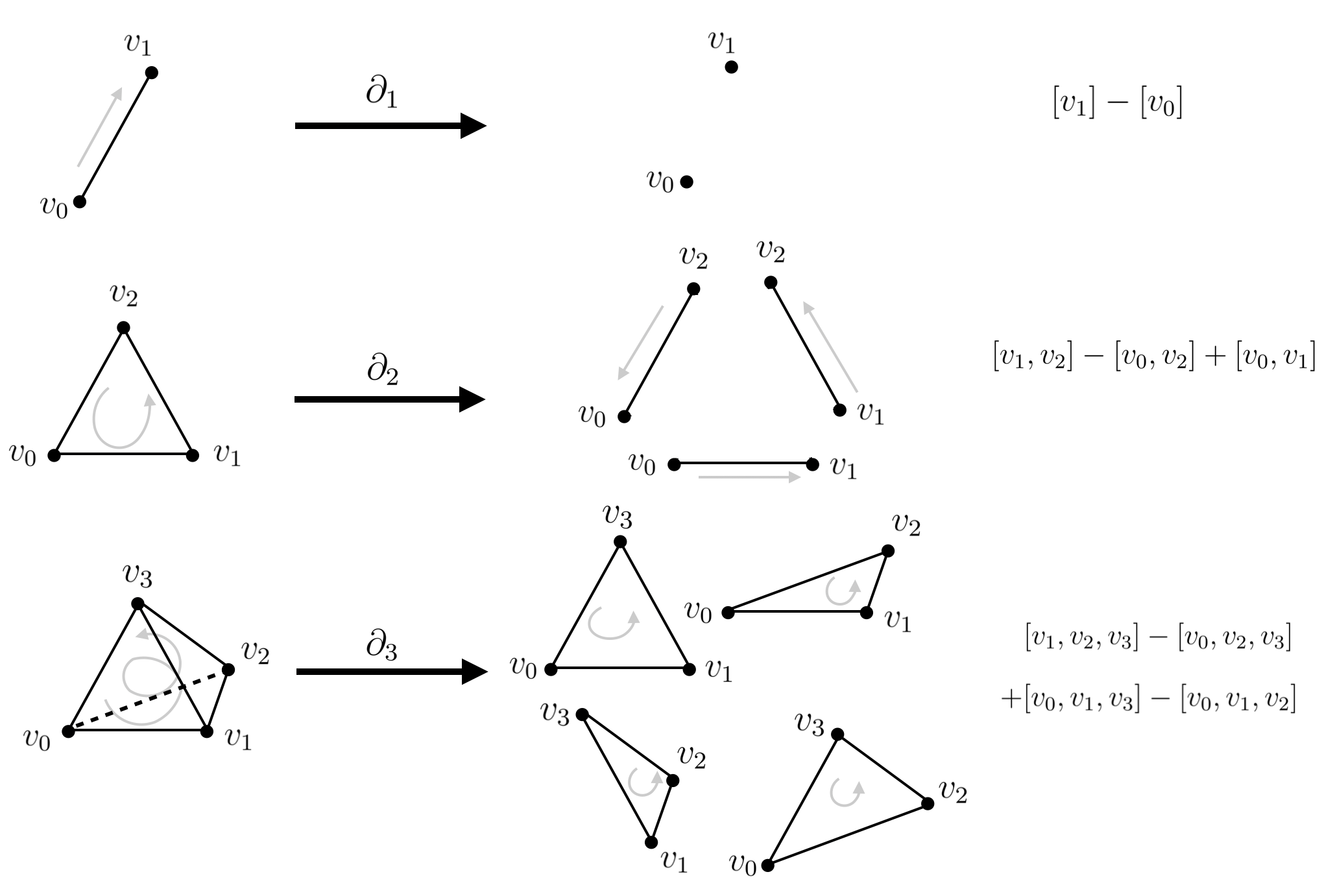}
    \caption{The linear map $\partial_n$ interpreted geometrically in terms of the boundary of a simplex and the induced orientations.}
    \label{fig:boundary}
\end{figure}

Let us now consider the linear map $\partial_n : C_n(S) \to C_{n-1}(S)$ defined on basis elements by
\begin{equation} \label{equation: boundary map}
\partial_n \big(\lbrack v_0, \dots, v_n \rbrack\big) = \sum_{i = 0}^n (-1)^i \lbrack v_0, \dots, \widehat{v}_i, \dots, v_n \rbrack,
\end{equation}
with $\widehat{v}_i$ denoting the absence of $v_i$ from the simplex. One can visualise $\partial_n$ geometrically in terms of the boundary of a basis element. Up to a sign, the basis elements appearing as summands on $\partial_n \big(\lbrack v_0, \dots, v_n \rbrack\big)$ are all simplices of dimension $n-1$ contained in $[v_0, \dots, v_n]$. Furthermore, \rev{the signs in this formula can be interpreted in terms of compatibility, or lack of it, among the orientations induced by the total orders on the vertices of the simplices involved. Figure~\ref{fig:boundary} presents some examples} 
The maps $\partial_n$ play a central role in \rev{a}lgebraic \rev{t}opology holding a significant amount of topological information. In particular, the difference between the dimension of the kernel of $\partial_n$ and the dimension of the image of $\partial_{n+1}$ is known as the \mbox{$n$-Betti} number of the simplicial complex, a powerful topological invariant generalising the Euler characteristic. For a systematic treatment of these ideas, please consult~\cite{hatcher2002topology}.

In this work we are not only interested in the topology of $S$, but also on the ``geometric'' structure encoded by 
weights assigned to its simplices. A \textit{weighted simplicial complex} is a pair $(S, w)$ where $S$ is a simplicial complex and $w = \{w_n : S_n \to \mathbb R_{ > 0}\}$ is a set of non-negative \textit{weight functions}.
In the following, a weighted standard simplex (see Example~\ref{example: standard simplex}) is referred to as a \textit{structural simplex}. 
 
\begin{example} \label{example: structural simplex from mutual information}
	We can build a structural simplex using a collection of random variables $X_0, \dots, X_N$ by following Example~\ref{example: standard simplex}, and defining the weights as follows. First, for each vertex $v_k$ in $\Delta^N$ we assign $w_0([v_k])=1$, and for each edge $[v_i,v_j]$ we assign the mutual information
	\begin{equation*}
	w_1([v_i, v_j]) = I(X_{v_i}, X_{v_j}).
	\end{equation*}
	Subsequently, for an $n$-simplex $[v_0, \dots, v_n]$ with $n > 1$, its weight is defined as the mean value of the mutual information of all pairs of random variables in $\{X_{v_0}, \dots, X_{v_n}\}$. \qed
\end{example}

Given a weighted simplicial complex $(S, w)$, one can encode the structural information provided by $w_n$ (where the subscript denotes the dimension) as an inner product on the vector space of $n$-signals on $S$. Explicitly, for any pair of $n$-simplices we have
\begin{equation}\label{eq:inner_prod}
\langle \sigma, \sigma^\prime \rangle_w =
\begin{cases}
w_n(\sigma) & \text{if } \sigma = \sigma^\prime, \\
0 & \text{otherwise}.
\end{cases}
\end{equation}
Importantly, this inner product allows us to introduce the adjoint operator of $\partial_n$, which is denoted by $\delta_n$. That is to say, the operator $\delta_n : C_n(S) \to C_{n+1}(S)$ is defined by the identity
\begin{equation*}
\langle \partial_{n+1} \alpha, \beta \rangle_w = \langle \alpha, \delta_{n} \beta \rangle_w
\end{equation*}
which holds for any $\alpha \in C_{n+1}(S)$ and $\beta \in C_n(S)$. Note that $\partial_n$ does not depend on the weights $w$, but $\delta_n$ does.  

\subsection{High-order Laplace operator and Fourier basis}
\label{sec:33}

The classical sine and cosine functions form the basis of the spectral representation of time-domain signals.
In effect, classical Fourier analysis guarantees that a large class of functions are expressible in terms of \rev{the} linear combination of these functions. The coefficients of this linear combination are the Fourier coefficients, which correspond to a geometric projection of the original function over this new basis. Interestingly, numerous signals of practical relevance are more compactly represented in the spectral domain, and hence Fourier analysis is often employed as a principled way to perform dimensionality-reduction. For a more extensive exposition of these ideas, we refer the reader to \cite{bracewell1986fourier}.

Importantly, sine and cosine functions are also eigenfunctions of the Laplace operator on the circle --- a one-dimensional manifold.
Put differently, a spectral representation is equivalent to a diagonalisation of the Laplace operator. Mathematicians have generalized the Laplace operator to higher-dimensional manifolds via the Laplace-de Rham operators. In this more general context, functions on a smooth manifold lie at the bottom of a sequence of higher-dimensional objects on which the corresponding Laplace-de Rham operator acts. The eigenvectors of these operators play a central role in modern geometry --- most notably through Hodge theory~\citep{hodge1989theory}. For an exposition of these ideas we refer the reader to \cite{morita2001geometry}.

A discrete analogue of the higher-dimensional objects on which the Laplace-de Rham operators act is provided by high-order signals defined over a simplicial complex. In this correspondence, the Laplace-de Rham operators are represented by the discrete Laplace operators first introduced in \cite{eckmann1944harmonische}.
See also \cite{horak2013spectra} where a weighted version of these is presented, and \cite{parzanchevski2017simplicial} where the connection to random walks is explored. We now introduce the particular version of the discrete Laplace-de Rham operators that we use.

\begin{definition} \label{definition: laplace operator}
	Let $(S, w)$ be a weighted simplicial complex. The $n$-\textit{Laplace operator} $\Delta_n : C_n(S) \to C_n(S)$ is defined by
	\begin{equation} \label{equation: discrete Laplace operator}
	\Delta_n = \partial_{n+1} \delta_n + \delta_{n-1} \partial_n.
	\end{equation}
\end{definition}

Notice that we are abusing notation by omitting $w$ from the notation, referencing the Laplace operators since, although $\partial_n$ does not depend on $w_n$, $\delta_n$ and therefore $\Delta_n$ do. For the interested reader, we remark that this Definition~\ref{definition: laplace operator} has a strong resemblance
\rev{to} the form
the Laplace-de Rham operator adopts in Rimmanian geometry. In effect, by denoting by $d$ the exterior derivative of differential forms and $d^\ast$ its adjoint, the Laplace-de Rham operator, in this context, can be expressed as $dd^\ast+d^\ast d$. Furthermore, the geometry defined by the Riemannian metric is reflected in the Laplace-de Rham operator \rev{by} the operator $d^\ast$ only. Consult \cite{morita2001geometry} for further details.

The well-known graph Laplacian, a central concept in spectral graph theory~\citep{chung1997spectral}, is equivalent to the $0$-Laplace operator defined above --- when a weighted graph is regarded as a weighted simplicial complex \rev{and the map $\partial_0$ is, as usual, set to be the $0$ map}.
Importantly, the $n$-Laplace operator is self-adjoint for any weighted simplicial complex $(S, w)$, i.e.
\begin{align*}
\big \langle \Delta_n(\alpha), \beta \big \rangle_w 
= & \ 
\big \langle \partial_{n+1} \delta_n (\alpha), \beta \big\rangle_w  
+ 
\big \langle \delta_{n-1} \partial_n (\alpha), \beta \big\rangle_w \\ 
= & \
\big \langle \alpha, \partial_{n+1} \delta_n (\beta) \big \rangle_w 
+
\big \langle \alpha, \delta_{n-1} \partial_n(\beta) \big \rangle_w 
\\ 
= & \
\big \langle \alpha, \Delta_{n}(\beta) \big \rangle_w
\end{align*}
for all $\alpha, \beta \in C_n(S)$. This implies that $\Delta_n(\alpha)$ is diagonalisable; that is to say, there exists a basis of $C_n(S)$ consisting of eigenvectors of $\Delta_n$.\footnote{This result is a more general form of the well-known fact that symmetric (i.e. self-adjoint) matrices are diagonalisable.} An $n$-\textit{Fourier basis} for $(S, w)$ is a maximal set of linearly independent orthonormal (with respect to $\langle\cdot,\cdot\rangle_w$) eigenvectors of $\Delta_n$. Please note that we speak of \textit{the} Fourier basis, with the understanding that there is a sign choice for each of its elements.

Finally, the \textit{hyperharmonic representation} of a high-order signal defined over a weighted simplicial complex is given by its change of bases from the canonical to the Fourier one. In the case of graphs, the use of this transformation as a dimensionality-reduction method was pioneered by \cite{belkin2003laplacian}, and has served as an early application of the field of graph signal processing --- see for example \cite{shuman2013emerging, sandryhaila2014discrete} or \cite{ortega2018graph} for an overview of this field. In recent years, some key ideas from graph signal processing have been adapted to hypergraphs, see for example \cite{barbarossa2020topological} and \cite{schaub2018flow}; however, the applicability of Fourier analysis as a compression tool of high-order signals is, to a large extend, still unexplored territory.

\section{Proposed workflow}
\label{sec:pipe}

\begin{figure}[h!]
    \centering
    \includegraphics[width=\textwidth]{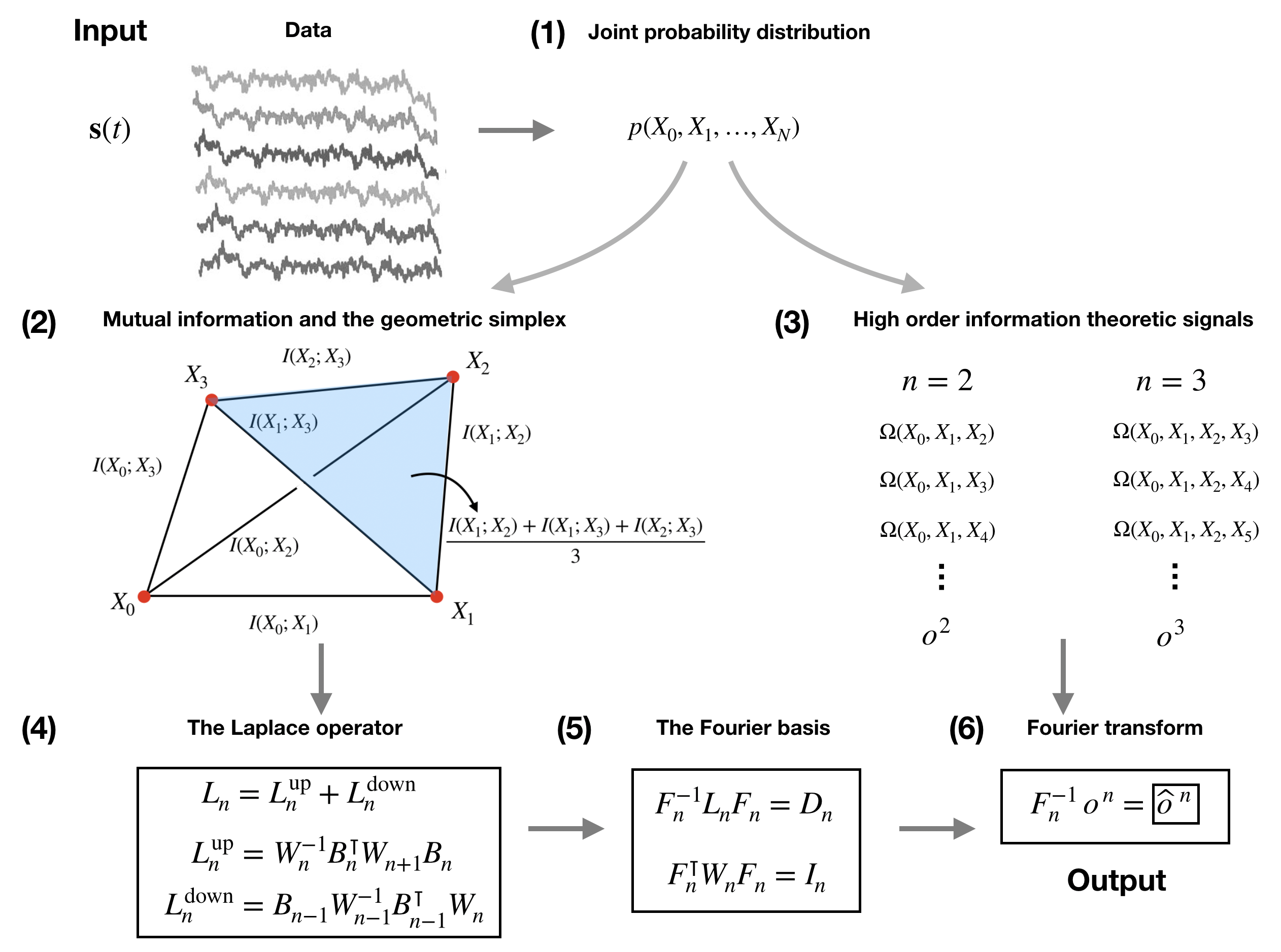}
    \caption{Workflow: (1) Estimate a joint probability from the input data (2) Build a structural simplex using averages of mutual information values. (3) Compute the \mbox{O-information} (and S-information) and store it as an \mbox{${N+1}\choose{n+1}$-dimensional vector} for every dimension $n$ (4) Compute the $n$-Laplace matrix using the boundary and weight matrices. (5) Diagonalise the Laplace matrices for each dimension. (6) Write the signal in the Fourier basis and return it as output.}
    \label{fig:wf}
\end{figure}

This section describes our proposed workflow, which capitalises \rev{on} the theoretical constructs elaborated in Sections~\ref{sec:hoitm} and \ref{sec:hyper}. The overall pipeline is illustrated in Figure~\ref{fig:wf}, being composed of six steps that are described in the following.

\subsection{Steps of the analysis}

The proposed pipeline consists \rev{of} six steps.

\begin{itemize}

\item[(1)] \textbf{From data to a distribution:} The starting point of the workflow is multivariate data, which typically takes the shape of a sequence of vectors of dimension $N+1$ (e.g. successive samples of time series of the form $\bi{s}(t) = \left(s_0(t), \dots, s_N(t)\right)$ at different values of $t\in\mathbb{R}$). Using these data, the first step is to construct a joint probability distribution for \rev{an} $(N+1)$-dimensional vector, denoted as $p(X_0, \dots, X_N)$. 

\item[(2)] \textbf{The structural simplex:} 
To build the structural simplex, the first step is to use $p$ to calculate the mutual information $I(X_i;X_j)$ for each pair of variables of the system. These values are then used to define the weights \rev{in} a pairwise network
\rev{with nodes corresponding to the variables $X_0, \dots, X_N$.}
Subsequently, a weighted $n$-simplex is \rev{built} for each $n = 2, \dots, N$ by taking the average value of all the edges that connect two elements within the simplex.

\item[(3)] \textbf{High-order signals:} For each subset of cardinality $n+1$ of $\{X_0, \dots, X_N\}$, one computes the O-information and S-information and arranges them as high-order signals following Example~\ref{example: omega and sigma as signals}. These signals are stored as \mbox{${N+1}\choose{n+1}$-dimensional vectors} $\bi{o}^n$ and $\bi{s}^n$, representing them in the canonical basis of simplices. Note that there is a standard order in the elements of this basis, which is given by the lexicographic principle (i.e. $[v_0, \dots, v_n] < [v'_0, \dots, v'_n]$ if $v_j < v'_j$ and $v_i = v'_i$ for all $i < j$).

\item[(4)] \textbf{The Laplace operator:} 
Using the weights of the structural simplex constructed in Step (2), one then builds the corresponding $n$-Laplace operator as in Definition~\ref{definition: laplace operator}. Concretely, we use the canonical basis to represent the $d =$ ${N+1}\choose{n+1}$ weights of order $n+1$ within a $d \times d$ diagonal matrix $W_n$, and represent the linear map $\partial_n$ as a matrix $B_n$ of the same dimensions (for an algorithmic description of the construction of $B_n$, see~\ref{appendix: contruct boundary matrix}). Then, the discrete $n$-Laplace operator is represented in the canonical basis as the matrix
\begin{equation}
L_n = L_n^{\text{up}} + L_n^{\text{down}},\label{eq: matrix laplace}
\end{equation}
with $L_n^{\text{up}}$ and $L_n^{\text{down}}$ given by
\begin{align*}
L_n^{\text{up}} & = W_{n}^{-1} B_n^\intercal \,W_{n+1} B_n, \\
L_n^{\text{down}} & = B_{n-1} W_{n-1}^{-1}\, B_{n-1}^\intercal W_{n},
\end{align*}
\rev{where} $B_n^\intercal$ is the transpose of $B_n$. To recapitulate, Eq.~\eqref{eq: matrix laplace} is the matrix representation, in the canonical basis, of Eq.~\eqref{equation: discrete Laplace operator} with respect to the inner product given by Eq.~\eqref{eq:inner_prod}.

\item[(5)] \textbf{The Fourier basis:} The $n$-Fourier basis of the structural simplex is constructed by choosing a \rev{maximal} linearly independent set of eigenvectors of the $n$-Laplace operator that are orthonormal with respect to the inner product $\langle \cdot,\cdot \rangle_w$ defined by the weights of the structural simplex (see Eq.~\eqref{eq:inner_prod}). Concretely, one needs to find a matrix $F_n$ such that
\begin{equation}
F_n L_n F_n^{-1} = D_n
\end{equation}
with $D_n$ a diagonal matrix, which also satisfies
\begin{equation*}
(F_n^{-1})^\intercal \, W_n  F_n^{-1} = I_n,
\end{equation*}
where $I_n$ is the identity matrix of dimension $N+1\choose{n+1}$.

\item[(6)] \textbf{Hyperharmonic representation:}
As a final step, one calculates the Fourier transform of the high-order signals calculated in Step (3), denoted by $\widehat{\bi{o}}^{\, n}$ and $\widehat{\bi{s}}^{\, n}$, as follows:
\begin{align*}
F_n \bi{o}^n & =  \widehat{\bi{o}}^{\, n}, \\ 
F_n \bi{s}^n & = \widehat{\bi{s}}^{\, n}.
\end{align*}
\end{itemize}

\subsection{Variations}\label{sec:variations}

This workflow has been designed with modularity and flexibility in mind, \rev{in order} to facilitate its adaptation to the needs of diverse applications. This section highlights possible variations \rev{in} the workflow that can better accommodate the specific needs of different scenarios.

First, note that our choice of $\Omega$ and $\Sigma$ in Step (3) as high-order information signals is based on their interpretability, and on their promising value for the analysis of complex systems. Nevertheless, other high-order measures (e.g. the ones described in \ref{subsec: other high order signals}) can also be analysed following the same pipeline. 

Additionally, Step (2) suggests \rev{building} the structural simplex by propagating the underlying mutual information from edges to higher-dimensional simplices via averages. However, one could use other constructions: e.g. use the maximum or minimum mutual information instead of the average.
Additionally, one could replace the mutual information with other non-negative metrics of similarity, including the total variation distance, the Wasserstein distance, or the absolute value of the Pearson correlation. Furthermore, one could also use non-negative high-order metrics (such as the TC or DTC, see Section~\ref{sec:Oinfo}) for building the structural simplex directly.

Finally, it is important to note that Step (1) is included mainly for pedagogical purposes, but is often \rev{omitted} in practice. In effect, most modern techniques to estimate information-theoretic quantities from data avoid \rev{building} an explicit joint distribution, as this introduces additional biases. For discrete data, we recommend considering Bayesian estimators such as the ones discussed in~\cite{archer2013bayesian}, and the software package DIT~\citep{james2018dit}. For continuous data, our \rev{preferred} choice \rev{is Kraskov} estimators~\citep{kraskov2004estimating} that can be implemented via JIDT~\citep{lizier2014jidt}. However, these estimators require substantial amounts of data, and a more flexible alternative \rev{is} provided by estimators based on  Gaussian Copulas~\citep{Ince2017}.

\section{Proof of concept: Haydn symphonies} \label{sec:example}

As an illustration of the proposed workflow, this section presents an analysis of the degree of dimensionality-reduction that can be attained by performing hyperharmonic analysis over the O-information and S-information signals calculated over a small dataset.
For this purpose, we use data from the music scores of Franz Joseph  Haydn (1732--1809), one of the most iconic figures of the Classical Period. We focus on Haydn's latter ``London symphonies'', which are typically divided into two groups: Symphonies 93--98, which were composed during the first visit of Haydn to London, and Symphonies 99--104, which were composed either in Vienna or London during Haydn's second visit~\citep{clark2005cambridge}. \rev{We employ this dataset in order to demonstrate the capabilities of our proposed framework on real data; however, the interpretation of these findings and their implications for music analysis is left to future work.}

\subsection{Method description}

Our analysis is based on electronic scores that are publicly available at
\url{http://kern.ccarh.org}, from where we extracted the scores of Symphonies No.~93--94 and No.~99--104.\footnote{The data for Symphonies 95--98 was not available.} All symphonies use the same basic instrumentation: flutes, oboes,
bassoons,
horns, trumpets, timpani, violins, viola, violoncello, and double bass;
only some of these symphonies employ clarinets, which were therefore not included in the analysis. To avoid \rev{duplicates},
our analyses consider only one
instrument of each kind, which left an arrangement of nine parts (\rev{violoncello}
and double bass
were assumed to be equivalent). 

The scores of the four movements of the selected symphonies were pre-processed in \texttt{Python~3.8.5} using the \texttt{Music21} package (\url{http://web.mit.edu/music21}). Each movement was transformed into nine coupled time series taking 13 possible values --- one for each note, plus one for the silence, using a small rhythmic duration\footnote{\rev{A small rhythm unit is chosen, so that each time unit contains only one pitch --- this analyses use a $1/48$-th of a quarter note as time unit.}} 
as \rev{a} time unit.
With these data, the joint distribution of the values for the nine-note chords was estimated using their
empirical frequency. Note that regularisation methods (such as Laplace smoothing) were not employed, as many configurations (e.g. highly dissonant chords) are never explored in the Classical repertoire. 

Our subsequent analyses were restricted to high-order signals depending on the joint probability of no more than six~instruments. Given the structure of the dataset (eight symphonies with four movements each), we computed the structural simplex using a distribution calculated using all \rev{of} the data. In contrast, individual high-order signals were calculated for each movement of each symphony by using only the corresponding data. 

To measure compressibility, we used the following metric defined for any basis. Consider a given $n$-dimensional high-order signal,
whose coefficients on the given basis are $\bm\alpha = \{\alpha_i\}$. Without loss of generality we assume that $\alpha_i^2 \geq \alpha^2_j$ if $i \leq j$. Then, we define the functions
\begin{equation}\label{eq:CEV}
   \text{EV}_{\alpha}(k) = \frac{ \alpha_k^2}{{\displaystyle \sum_{i} \alpha_i^2}}
   \qquad \text{ and } \qquad
   \text{CEV}_{\alpha}(k) = \sum_{1\leq i \leq k} \text{EV}_{\alpha}(i),
\end{equation}
with $\text{EV}_\alpha(k)$ being the \textit{(normalised) explained variance} by the $k$-th strongest component, and $\text{CEV}_\alpha(k)$ the \textit{cumulative explained variance} by the $k$ strongest components. This definition is motivated by Parseval's theorem, which guarantees that the sum of the square of the Fourier coefficients is equal to the variance of the signal; hence, the $\text{CEV}(k)$ of the Fourier transform of a signal corresponds to the percentage of the variance that is accounted by the $k$ strongest components.

\subsection{Results}
\label{sec:results}

\begin{figure}[t!]
    \centering
    \includegraphics[width=\textwidth]{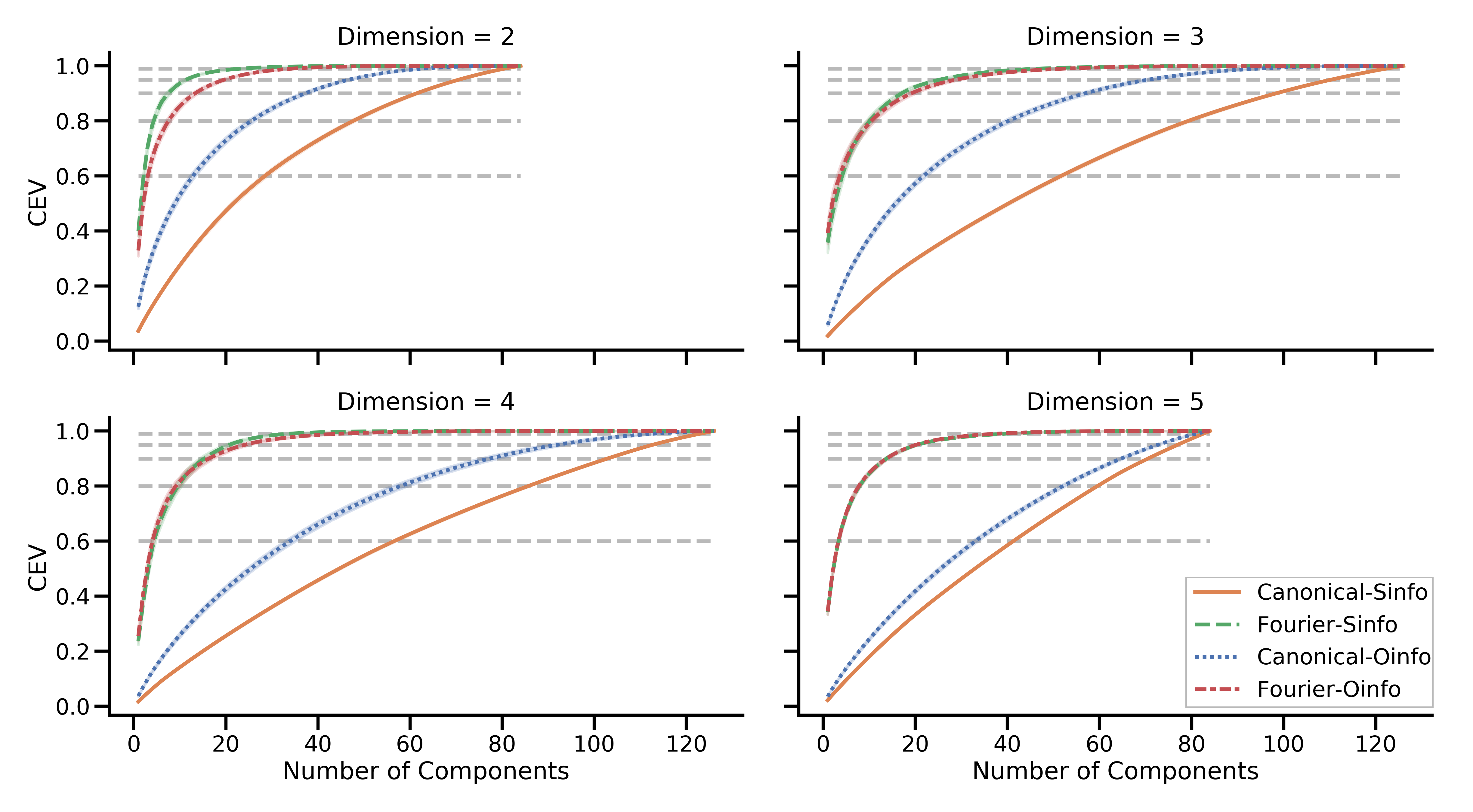}
    \caption{Comparison \rev{(via the average and 95\% confidence interval)} of the cumulative explained variance (as defined by Eq.~\eqref{eq:CEV}) of high-order signals in their canonical and hyperharmonic representation.
    The dimensionality-reduction capabilities of hyperharmonic analysis is substantial, and stays consistent across dimensions.}
    \label{fig:CEV}
\end{figure}

The hyperharmonic representation of the high-order signals considered ($\Omega$ and $\Sigma$) was found to be substantially more concentrated than their corresponding canonical representations. Figure~\ref{fig:CEV} illustrates the curves of \rev{the}  cumulative explained variance for various dimensions, and Table~\ref{table} presents the number of components required to fulfil various reconstruction levels. In particular, it was found that a small number of components suffices to account for most of the variance observed in the hyperharmonic representations of $\Omega$ and $\Sigma$.  \rev{Moreover, this effect is found to be stable across dimensions.}

\begin{table}[ht!]
\centering
\begin{tabular}{ |p{1.6cm}|p{1cm}|p{1cm}|p{1cm}| p{1,2cm}| p{1,3cm}| p{1,3cm}|}
 \hline \hline
 \multicolumn{7}{|c|}{\textit{Cumulative Explained Variance}} \\
 \hline
\textit{Signal} & \textit{Dim}  
& 60\% & 80\% & 90\% & 95\% & 99\%\\
 \hline\hline
\multirow{3}{*}{\textbf{O-info}}
& 2 &    3\,/\,13   &   5\,/\,26   &  9\,/\,38    &  14\,/\,47     &  28\,/\,63    \\

& 3 &     4\,/\,22   &   8\,/\,41   &  13\,/\,57   &  20\,/\,71    &  37\,/\,95   \\

& 4 &     1\,/\,35   &   3\,/\,58   &  6\,/\,78    &  12\,/\,93     &  29\,/\,113  \\

& 5 &     2\,/\,34   &   4\,/\,53   &  8\,/\,65    &  12\,/\,73     &  26\,/\,82  \\
   \hline
\multirow{3}{*}{\textbf{S-info}}
& 2 &    2\,/\,29   &   4\,/\,48   &  6\,/\,62    &  8\,/\,71     &  19\,/\,81    \\

&3 &    3\,/\,52   &   7\,/\,80   &  11\,/\,99   &  14\,/\,111   &  27\,/\,123   \\

& 4 &    1\,/\,57   &   2\,/\,86   &  3\,/\,103   &  9\,/\,113    &  24\,/\,123   \\

& 5 &    2\,/\,42   &   4\,/\,60   &  8\,/\,71   &  13\,/\,78    &  26\,/\,83   \\
 \hline\hline
\end{tabular}
\caption{Number of components needed to recover a given percentage of the cumulative explained variance for the considered signals, either in the Fourier basis or in the canonical basis.}
\label{table}
\end{table}%

To verify the unique properties of the Laplace operators, an additional control was run \rev{where} the same signals were transformed according to randomly-generated bases. Interestingly, these representations do not exhibit the degree of dimensionality-reduction shown by the hyperharmonic representations~(see~\ref{app:control}).

Finally, our results also show that the canonical representation of the S-information is much less concentrated than the canonical representation of the O-information. This dissimilarity suggests that the O-information is capturing a more specific signature of the different modalities of interdependency that exist across the orchestra. Moreover, the fact that the dissimilarity is not seen in their hyperharmonic representation further suggests that both signals may actually carry an equivalent amount of information, which happens to be localised differently over the canonical basis.

\section{Conclusion}\label{sec:conc}

Complex systems are characterised by having multiple levels of organisation, which makes network analyses focused on pairwise interactions unable to give a full account of their \rev{properties}. 
Phenomena beyond pairwise interactions can be effectively captured by high-order information-theoretic metrics; however, their applicability and interpretability is limited due to the fast growth of their cardinality as a function of the system's size. \rev{Here w}e propose to represent these high-order metrics using the hyperharmonic modes of a geometrical representation of structural properties of the system. This provides a principled approach to \rev{constructing} 
low-dimensional representations of high-order signals, which 
\rev{can retain} 
most of the informational content.

Our proposed approach is widely applicable, and promises to enable a range of future explorations. \rev{It is to be noted that while the presented framework focuses on synchronous interdependencies, diachronous relationships could also be considered via recent extensions of the O-information introduced by~\cite{stramaglia2021quantifying}. Fully developing such approaches, and clarifying how this could be combined with ongoing work in causal discovery (e.g. Refs~\citep{runge2019detecting,liu2021quantifying}), is a promissing avenue for future investigation.
Additionally, future work might explore computationally efficient implementations of our proposed algorithm (e.g. by only partially diagonalising the Laplace operator matrix), which may allow a more graceful scaling of the computational complexity of of the proposed algorithms with respect to the system's size.}

It is our hope that this incursion into the intersection of multivariate information theory and combinatorial topology will motivate further developments in this fertile area of research.

\section*{Data Availability}
\rev{The data that support the findings of this study are openly available.}

\ack
The authors thank Pedro A. M. Mediano and Umberto Lupo for insightful discussions and valuable suggestions. The authors also thank Kathryn Hess and Yike Guo for supporting this research. A.M-M. acknowledges financial support from Innosuisse grant \mbox{32875.1 IP-ICT - 1}, and the hospitality of the Max Plank Institute for Mathematics. R.C. acknowledges financial support from Fondecyt Iniciaci\'{o}n 2018 Proyecto 11181072. F.R. acknowledges the support of the Ad Astra Chandaria Foundation. 

\appendix	

\section{Additional experiments}

\subsection{Control employing randomly generated bases}
\label{app:control}

\begin{figure}[h]
    \centering
    \includegraphics[width=\textwidth]{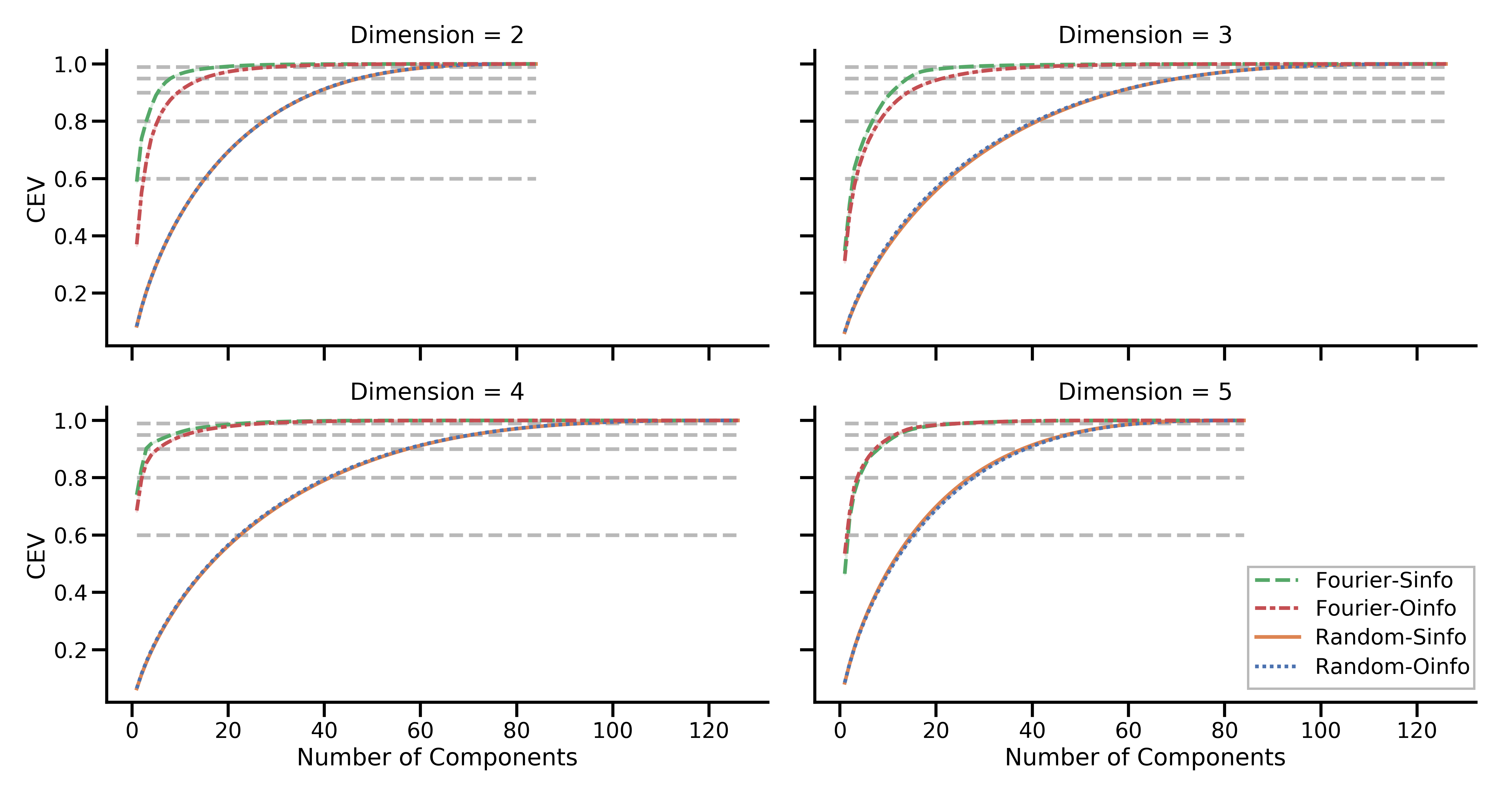}
    \caption{Comparison \rev{(via the average and 95\% confidence interval)} of the cumulative explained variance (as defined by Eq.~\eqref{eq:CEV}) of high-order signals in the hyperharmonic representation and the average CEV of these signals over randomly generated bases. The dimensionality-reduction capabilities of hyperharmonic analysis is substantial, and stays consistent across dimensions.}
    \label{fig:CEV-Random}
\end{figure}

Here we review the additional control that confirms that the results presented in Section~\ref{sec:results} are not due to a limitation of the canonical bases, but are due to the special advantages of the Laplace operators. Specifically, we calculated the average CEV over $80$ randomly generated bases for the same high-order signals considered in Section~\ref{sec:results}, and then compared them to the CEV obtained via the hyperharmonic representation. Our results, shown in Figure~\ref{fig:CEV-Random}, reveal that the explained variance attainable via randomly generated bases is substantially and consistently lower than the one associated to the Fourier bases. This provides additional evidence on the favorable properties of the hyperharmonic representation of the considered signals.

\subsection{\rev{Control employing synthetic signals}}
\label{app:signals}

\rev{For deepening our intuition about the properties of our method, 
we analised the dimensionality-reduction capabilities of the 
pipeline proposed in Section~\ref{sec:pipe} on synthetically-generated data. For this purpose we consider data built with internal correlations
having a controlled number of degrees of freedom, which are related to the dimensionality of their effective phase space.
We set out to evaluate how this type of ``diversity" on
the input data impacted the performance of our dimensionality-reduction method;
we conjectured that a system with less diverse input data would be more
compressible since, intuitively, it has less information distributed over the same number of simplices.}

\rev{To test our conjecture, we considered a scenario of
nine time series (i.e. $N=8$) with a joint probability distribution 
$p(X_0,\dots, X_8)$ given by a multivariate Normal distribution
with zero mean and a random covariance matrix $C_r$ of rank $r \in \{1,\dots,9\}$.
The rank of $C_r$ was used as a proxy for the diversity of the input, as a low-rank
distribution is constrained to a low-dimensional subspace of $\mathbb{R}^9$. 
For building $C_r$ we employed the following procedure:
(i) generate a random positive-definite matrix $A$ by generating a $9\times 9$ random matrix $M$ with I.I.D.
components following a standard Normal distribution, and then calculating $A= M M^T$,
(ii) calculate the singular value decomposition of $A$ as $A = V \Sigma V^T$,
(iii) build $\Sigma_r$ by turning to zero the $9-r$ smaller eigenvalues, and
(iv) compute $C_r = V \Sigma_r V^T$.
After building $C_r$, we sampled $10^4$ values of $(X_0,\dots,X_8)$ following
the corresponding distribution, which
were then used as input for the procedure described in Sect~\ref{sec:pipe}.
For reconstructing an empirical joint distributions, we employed the method of Gaussian copulas~\citep{Ince2017}.
Finally, for accounting the variability of $C_r$, we repeated this process $50$ times, and then reported the average CEV (as given in Eq.~\eqref{eq:CEV}) over all realisations.}

\begin{figure}[h]
    \centering
    \includegraphics[width=0.9\textwidth]{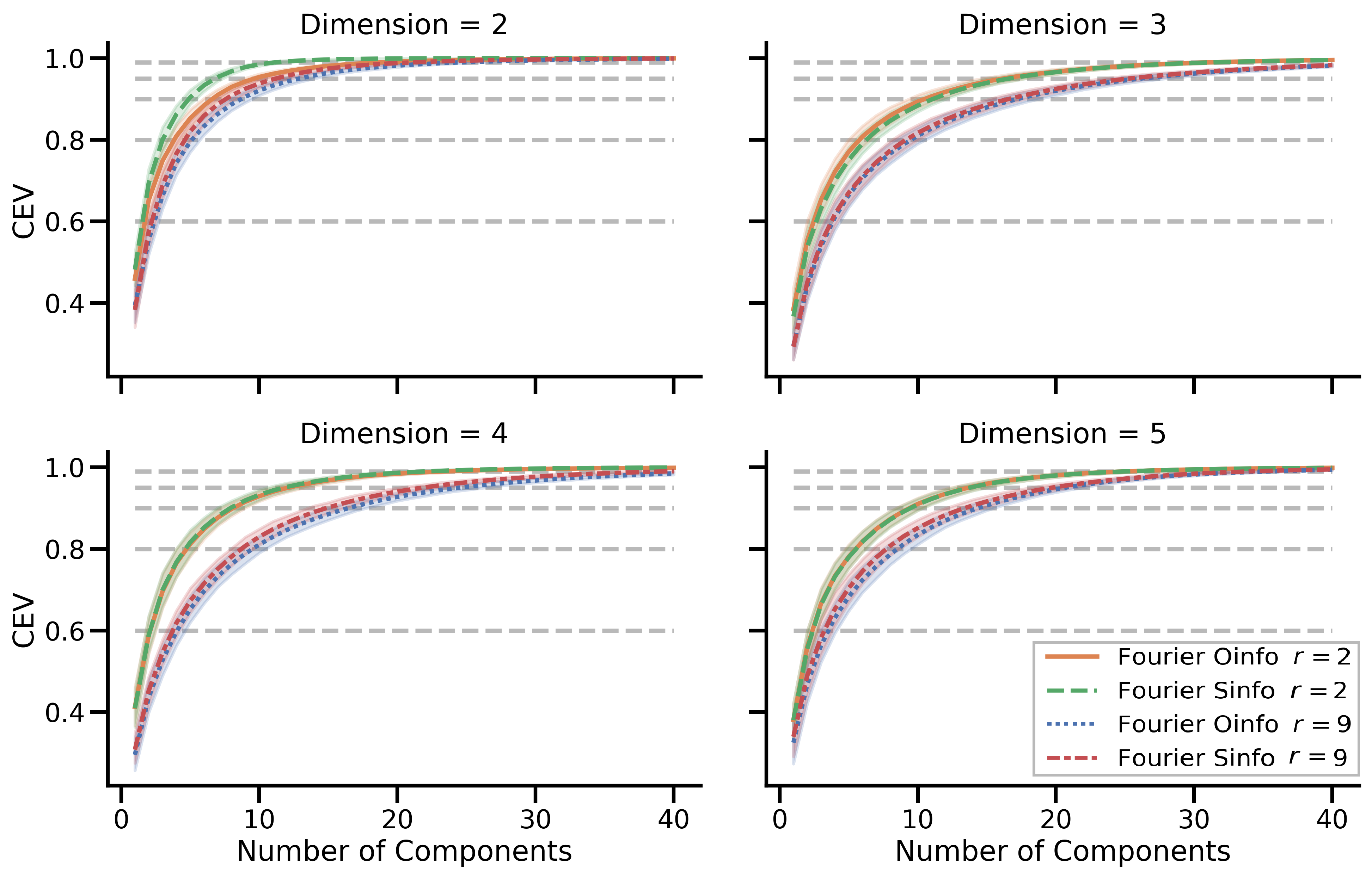}
    \caption{Comparison \rev{(via the average and 95\% confidence interval)} of the cumulative explained variance across dimensions and values $r=2$ and $r=9$ of the rank of the correlation matrix used to generate the input data. The dimensionality-reduction capabilities of our hyperharmonic method are more prominent for those with small rank.}
    \label{fig:CEV_appendix}
\end{figure}

\rev{Our results, shown in Figure \ref{fig:CEV_appendix}, confirmed the stated hypothesis, showing that systems with 
lower rank are substantially more compressible than those with full rank. 
As expected, this effect is more noticeable at dimensions $3$
and $4$, where the number of simplices is larger.
}

\section{Construction of $B_n$}
\label{appendix: contruct boundary matrix}

The matrices $B_n$ correspond to the representation of the linear maps $\partial_n \colon C_n(\Delta^N) \to C_{n-1}(\Delta^N)$ on the canonical basis (see Section~\ref{sec:32} and Figure~\ref{fig:boundary}). Algorithmically, one can use Eq.~\eqref{equation: boundary map} to determine the value of each column \rev{by} inserting either a $+1$ or a $-1$ to the entries corresponding to simplices in its boundary --- as depicted in Figure \ref{fig:boundary} --- and setting the other entries to~0. 
To illustrate the procedure, let us present the four matrices that correspond to $N=3$:
\begin{align*}
B_0 &= \begin{bmatrix}	0 & 0 & 0 & 0
	\end{bmatrix}, &
B_1 &= \begin{bmatrix}
	-1 & -1 & -1 & \phantom{+}0 & \phantom{+}0 & \phantom{+}0\\
	+1 & \phantom{+}0 & \phantom{+}0 & -1 & -1 & \phantom{+}0\\
	\phantom{+}0 & +1 & \phantom{+}0 & +1 & \phantom{+}0 & -1\\
	\phantom{+}0 & \phantom{+}0 & +1 & \phantom{+}0 & +1 & +1\\
	\end{bmatrix}, \\
B_2 &= \begin{bmatrix} +1 & +1 & \phantom{+}0 & \phantom{+}0\\
	-1 & \phantom{+}0 & +1 & \phantom{+}0\\
	\phantom{+}0 & -1 & -1 & \phantom{+}0\\
	+1 & \phantom{+}0 & \phantom{+}0 & +1\\
	\phantom{+}0 & +1 & \phantom{+}0 & -1\\
	\phantom{+}0 & \phantom{+}0 & +1 & +1\\
	\end{bmatrix}, &
B_3 &= \begin{bmatrix}
	-1\\
	+1\\
	-1\\
	+1\\
	\end{bmatrix}.
\end{align*}

\bibliography{main}
\end{document}